 \documentclass[final,1p,times,twocolumn]{elsarticle}

 \usepackage{graphics}

\usepackage{amsmath}
\usepackage{xcolor}
\usepackage{amssymb}
\usepackage{amsfonts}
\usepackage{fancybox}
\usepackage{color}
\usepackage{bm}
\usepackage{cancel}

\usepackage{graphicx}
\usepackage{epstopdf}

\usepackage{amssymb}
\usepackage{graphicx}
\usepackage{wrapfig}
\usepackage{caption}
\usepackage{subcaption}
\usepackage{amsmath}
\usepackage{graphicx,color}
\usepackage{bbm}
\usepackage{amsmath}
\usepackage{dsfont}
\usepackage{multirow}
\usepackage{graphicx, color}
\usepackage{amsmath,amsfonts,amssymb}
\DeclareMathOperator*{\argmin}{arg\,min}

\usepackage{amsthm}

\journal{Journal of Computational Physics}

\begin{document}

\begin{frontmatter}



\title{A compatible high-order meshless method for the Stokes equations with applications to suspension flows}


\author[label1]{Nathaniel Trask}
\author[label1]{Martin Maxey}
\author[label2]{Xiaozhe Hu}

\address[label1]{Division of Applied Mathematics, 182 George St., Providence, RI 02912}
\address[label2]{Department of Mathematics, Tufts University, 503 Boston Ave., Medford, MA 02155}

\begin{abstract}
A stable numerical solution of the steady Stokes problem requires compatibility between the choice of velocity and pressure approximation that has traditionally proven problematic for meshless methods. In this work, we present a discretization that couples a staggered scheme for pressure approximation with a divergence-free velocity reconstruction to obtain an adaptive, high-order, finite difference-like discretization that can be efficiently solved with conventional algebraic multigrid techniques. We use analytic benchmarks to demonstrate equal-order convergence for both velocity and pressure when solving problems with curvilinear geometries. In order to study problems in dense suspensions, we couple the solution for the flow to the equations of motion for freely suspended particles in an implicit monolithic scheme. The combination of high-order accuracy with fully-implicit schemes allows the accurate resolution of stiff lubrication forces directly from the solution of the Stokes problem without the need to introduce sub-grid lubrication models.
\end{abstract}


\end{frontmatter}


\section{Introduction}\label{sec:Introduction}

In this work we consider the steady Stokes problem
\begin{equation}
  \begin{cases}
    -\nu \nabla^2 \mathbf{u} + \nabla p = \mathbf{f} & \text{if } x \in \Omega\\
    \nabla \cdot \mathbf{u} = 0  & \text{if } x \in \Omega\\
    \mathbf{u} = \mathbf{w}  & \text{if } x \in \partial \Omega\\
  \end{cases}
  \label{StokesEqn0}
\end{equation}
where $\Omega \subset \mathbb{R}^d$ has a piecewise continuous boundary $\partial \Omega$, $\mathbf{u}$ and $p$ are velocity and pressure, $\nu$ is the kinematic viscosity, and $\mathbf{f}$ and $\mathbf{w}$ are given data and Dirichlet conditions. Meshless discretizations have long promised an attractive framework for discretizing fluid mechanics problems, due to their ability to handle easily large boundary deformations, to track explicitly interfacial phenomena in flows, and to discretize Lagrangian or semi-Lagrangian schemes. In this work we present a new strong form meshless discretization of the Stokes problem that is able to achieve equal high-order convergence for both velocity and pressure while maintaining a sparse discretization.

When solving the Stokes system, all numerical methods require some form of compatibility to control spurious solutions in the pressure field. For example, in finite element methods this amounts to an inf-sup condition between the velocity and pressure spaces, typically requiring a lower order approximation of pressure than velocity \cite{braess2007finite}. Finite difference and finite volume methods require a staggered arrangement of velocity and pressure variables or a separate mesh for velocity and pressure\cite{harlow1965numerical,mckee2008mac,ferziger2012computational,chou1998mixed}. Recently the discrete exterior calculus has unified the concept of compatibility for these different methods by considering the discretization of exterior derivatives when the mesh is interpreted as a simplicial complex\cite{arnold2007compatible,arnold2006finite,arnold2010finite,bochev2006principles,perot2011discrete}. Meshless methods, on the other hand, have no intuitive notion of a mesh and as a result cannot make use of this framework. Some methods that act on point sets are able to achieve a sense of compatibility by reintroducing some notion of a mesh. For analysis purposes, staggered finite differences on a Cartesian grid may be interpreted as a discretization of the staggered finite volume method on a rectilinear grid. In the mimetic finite difference method, Voronoi cells are introduced to obtain a simplicial complex and define discretely adjoint discretizations\cite{lipkikov2014mimetic}. 

Smoothed particle hydrodynamics (SPH) is one of the oldest meshless methods, and due to a notion of conservation has been successfully applied to a range of applications despite a lack of even zeroth-order consistency in some formulations\cite{lucy1977numerical,gingold1977smoothed,monaghan1983shock,monaghan2012smoothed,price2008modelling}. Somewhat recently the method has been used to discretize projection methods for the Stokes and Navier-Stokes equations\cite{cummins1999sph,xu2009accuracy}. Several authors have sought a consistent high-order discretization by applying either a least squares reconstruction or kernel correction to restore polynomial consistency; see for example moving least squares (MLS)\cite{dilts1999moving}, reproducing kernel particle method\cite{liu1995reproducing}, finite point method \cite{onate1998mesh,onate1996stabilized,onate1996finite}, and generalized finite differences\cite{gavete2003improvements}.  For these approaches, particles carry a velocity and pressure in the same location and therefore resemble unstable collocated finite difference methods, albeit defined on a scattered pointset. Successful discretizations must therefore either explicitly introduce stabilization or use variants of Chorin's first order projection scheme, which has been shown in a finite element context to implicitly introduce similar stabilizing terms\cite{guermond2006overview,liu2007stability}. The authors' own experience using a higher-order projection scheme lacking stabilization has revealed spurious pressure solutions at very high and very low Reynolds numbers \cite{traskCMAME}. The current work was motivated by an attempt to use these projection methods to study flows with suspended rigid particles. For these problems high accuracy is required to resolve singular lubrication effects and, as projection methods contaminate the convergence rate of pressure, a stable solution of the full Stokes system is required.

Without introducing stabilizing terms, few meshless options are available for the Stokes problem. Some methods have found success by introducing specialized point sets, using two staggered point sets for pressure and velocity or by using separate sets of data sites and collocation points\cite{kim2003point}. Within the framework of radial basis functions, Wendland obtained high-order approximation of the Stokes system by generating analytically divergence-free approximation spaces for velocity using matrix-valued kernels \cite{wendland2009divergence,narcowich1994generalized}. While this approach is highly accurate, it requires the inversion of a dense interpolation matrix and therefore extensions to large scale problems requires multi-level approaches to obtain a scalable solution \cite{farrell2013rbf}. Recent work by Fuselier \cite{fuselier2015high} has investigated the use of similar kernels to generate Leray projections that could be used to decouple the Stokes problem without introducing splitting errors. But these too require the solution of dense matrices. We will demonstrate that our approach is able to obtain similarly accurate results with a sparse approximation that can be solved using standard algebraic multigrid (AMG) techniques.

Dense suspension flow problems provide an ideal motivating application for such a meshless approach\cite{guazzelli2011physical,leal1980particle}. As colloidal particles are advected by flow, the geometry of the domain undergoes severe deformation that makes the generation of a colloid conforming mesh at every timestep very challenging. At the same time, as particles approach one another singular lubrication forces lead to singularities in the pressure field. The length scale of this so-called lubrication gap between particles can be less than $1\%$ of the particle diameter and presents a difficult multi-resolution problem for a variety of numerical methods. Traditional Eulerian approaches used by the suspension flow community can be broadly categorized as either: using a diffuse representation of the colloid boundaries\cite{peskin2002immersed,lai2000immersed,luo2009smoothed,glowinski2001fictitious,patankar2000new}, using Green's functions for spherical colloids to build representations of particles as point forces\cite{brady1988stokesian,lomholt2003force,yeo2010simulation}, using low-order meshless methods with artificial compressibility to avoid solving the coupled Stokes problem \cite{aidun2010lattice,nguyen2002lubrication,bian2014splitting,hashemi2012modified,hashemi2011sph}, or using ALE-type methods with a mesh conforming to colloid surfaces but requiring expensive remeshing and small timesteps\cite{hu2001direct,hu1996direct}. Many of these approaches are effective for this problem but introduce restrictions: diffuse interface approaches are usually at best second-order approximations that require sufficient resolution across the diffuse interface, while Green's function-type methods are restricted to simplified geometries. Many of these techniques rely on sub-grid scale lubrication models to avoid resolving the possibly narrow gap between particles and short-range repulsive potentials are used to represent unresolved colloidal forces or effects of rough surface contacts. We will demonstrate that our new approach is able to resolve this challenging multiscale problem directly by using a combination of high-order pressure approximation and adaptivity. In application our scheme would likely also incorporate lubrication corrections to obtain cost-effective results. We omit such corrections here to demonstrate the high-order of accuracy and sharp representation of boundaries that our scheme allows.

This paper is organized as follows. We first define a model problem in Section \ref{sec:colloidEquations} demonstrating how the Stokes problem can be coupled together with the equations of motion for suspended colloidal particles to obtain a fully implicit solver. This monolithic approach simultaneously couples the colloid translational and rotational velocities together with the fluid velocity and pressure. In Section \ref{sec:implementation} we present the details of the new discretization. First, in \ref{subsec:discretize} we describe our new scheme that combines a novel divergence-free MLS reconstruction of velocity with a staggered reconstruction for the pressure terms to obtain finite difference-like stencils at each particle. Second, in \ref{subsec:multiresolution} we describe a scheme to introduce adaptivity into the discretization. Thirdly, in \ref{subsec:forceCalc} we present details of how to evaluate stresses at the colloid surface.  Subsequently, in \ref{subsec:solvers} we provide a description of how AMG can be used to efficiently solve the resulting system. We then systematically validate each component of the scheme. In Section \ref{sec:stokesval} we present benchmarks demonstrating that when solving Equation \ref{StokesEqn0} for fixed curvilinear geometries, our approach is able to obtain equal high-order convergence for both velocity and pressure, and that the viscous stress approximation in \ref{subsec:forceCalc} is accurate for well-resolved lubrication forces. In Section \ref{sec:results} we demonstrate that the monolithic solver is accurate for channel flow and shear flow cases for which semi-analytic benchmarks for colloid drift velocity and trajectories are available.

\section{Colloid model problem}\label{sec:colloidEquations}

We consider in this work the motion of small, freely suspended particles under conditions for which Brownian motion is negligible. To study such suspension flows, we provide the following additonal assumptions to Equation \ref{StokesEqn0}. The domain $\Omega$ contains $N_c$ colloidal particles each with boundary $\partial {\Omega}_i$. The boundary of the domain can therefore be partitioned as $\partial \Omega = \partial \Omega_D \cup \left( \cup_{i} \partial \Omega_i \right)$, where a flow $\mathbf{w}$ is prescribed along the portion of the boundary $\partial \Omega_D$. Each colloid has position $\mathbf{X}_i$, orientation $\mathbf{\Theta}_i$, and undergoes a solid body rotation characterized by translational and rotational velocities $\dot{\mathbf{X}}_i$ and $\dot{\mathbf{\Theta_i}}$. Given $\mathbf{w}$ specified at the boundary, the flow within $\Omega$ is completely specified by $\left\{\mathbf{X}_i,\mathbf{\Theta}_i\right\}_{i=1,\dots,N_c}$ and, assuming steady Stokes flow, the fluid state is coupled to the colloid motion via

\begin{equation}
  \begin{cases}
    -\nu \nabla^2 \mathbf{u} + \nabla p = f & \text{if } x \in \Omega\\
    \nabla \cdot \mathbf{u} = 0  & \text{if } x \in \Omega\\
    \mathbf{u} = \mathbf{w}  & \text{if } x \in \partial \Omega_D\\
    \mathbf{u} = \dot{\mathbf{X}}_i + \dot{\mathbf{\Theta}}_i \times (x - \mathbf{X}_i)   & \text{if } x \in \partial \Omega_i, \text{for all }i\\
  \end{cases}
  \label{StokesEqn1}
\end{equation}

If the velocity is chosen from an appropriate space of divergence-free vector fields, this is equivalent to solving the following pair of equations coupled only through an inhomogeneous Neumann pressure boundary condition.

\begin{equation}
  \begin{cases}
    \nu \nabla \times \nabla \times \mathbf{u} + \nabla p = f & \text{if } x \in \Omega\\
    \mathbf{u} = \mathbf{w}  & \text{if } x \in \partial \Omega_D\\
    \mathbf{u} = \dot{\mathbf{X}}_i + \dot{\mathbf{\Theta}}_i \times (x - \mathbf{X}_i)  & \text{if } x \in \partial \Omega_i, \text{for all }i\\
    \end{cases}
  \label{StokesEqn2a}
\end{equation}

\begin{equation}
\begin{cases}
    \nabla^2 p = \nabla \cdot \mathbf{f}                                                       & \text{if } x \in \Omega\\
    \partial_n p + \nu \hat{n} \cdot \nabla \times \nabla \times \mathbf{u}= \hat{n} \cdot f   & \text{if } x \in \partial \Omega\\
  \end{cases}
\label{StokesEqn2b}
\end{equation}
following the vector identity $\nabla^2 \mathbf{u} = -\nabla \times \nabla \times \mathbf{u} = \nabla \nabla \cdot \mathbf{u}$.

By recasting the Stokes operator in this equivalent form, the saddle-point structure of the Stokes operator is avoided, and after discretization we will obtain instead a system with elliptic matrices along the diagonal blocks that is more amenable to standard preconditioning techniques. The use of a curl-curl viscous operator relies on the introduction of a divergence-free velocity reconstruction that will be introduced in the following section.

The translation and angular velocity of each particle may be computed explicitly from the equations of motion

\begin{equation}
\begin{cases}
    m_i \ddot{\mathbf{X}_i} = \int_{\partial \Omega_i} \mathbf{\sigma}_V \cdot d\mathbf{A} & \text{for all }i \\
    I_i \ddot{\mathbf{\Theta}_i} = \int_{\partial \Omega_i}  \left( x - \mathbf{X}_i\right) \times \mathbf{\sigma}_V \cdot d\mathbf{A}  & \text{for all }i. \\
    \end{cases}
\label{forcetorque}
\end{equation}
Here $m_i$ and $I_i$ are the mass and moment of inertia, respectively, of particle $i$, and $\mathbf{\sigma}_V = -p \mathbf{I} + \frac{\nu}{2}\left( \nabla \mathbf{u} + \nabla \mathbf{u}^\intercal\right)$ is the viscous stress tensor. These velocities would then provide Dirichlet boundary conditions for \ref{StokesEqn2a} in a separated system for the flow and particles. If the inertia of the particles is negligible, for example for neutrally buoyant particles in a viscous Stokes flow, and they are force-free and torque-free then \ref{forcetorque} become constraints of zero hydrodynamic force and torque. This assumption is made here to match cases for which analytical solutions are available.  

For the purposes of this work, we will first study benchmarks with specified colloid translations and rotations to demonstrate the high order convergence of the discretization. We will then solve Equations \ref{StokesEqn2a},\ref{StokesEqn2b}, and \ref{forcetorque} concurrently as a monolithic system to study the dynamics of colloidal suspensions under the assumption that the flow is completely in equilibrium. Therefore, given a set of colloid states $(\mathbf{X}_i,\mathbf{\Theta}_i)_{i=1,\dots,N_C}$, the coupled solution of these equations provides the complete dynamics of the problem. If we write the set of colloid degrees of freedom as $\vec{\mathbf{X}} = \left\{ \mathbf{X}_1 , \dots , \mathbf{X}_{N_c},\mathbf{\Theta}_1,\dots,\mathbf{\Theta}_{N_c}\right\}$, we can concisely describe the dynamics of the colloid motion as
\begin{equation}
\dot{\vec{\mathbf{X}}} = F(\vec{\mathbf{X}},t)
\label{colloidEvolution1}
\end{equation}
where $F$ denotes the solution of the given Stokes system with boundary prescribed by $\vec{\mathbf{X}}$ and driven by a flow $\mathbf{w}$. Equation \ref{colloidEvolution1} can then be integrated in time to evolve the particle locations and calculate particle trajectories. To perform this integration, we apply a second order Adams-Bashforth/Adams-Moulton predictor-corrector scheme.

\begin{align}
\frac{ \vec{\mathbf{X}}^{*}  - \vec{\mathbf{X}}^{n} }{\Delta t} &= \frac{ 3 F(\vec{\mathbf{X}}^n,t^n) - 3 F(\vec{\mathbf{X}}^{n-1},t^{n-1})}{2}\\
\frac{ \vec{\mathbf{X}}^{n+1} - \vec{\mathbf{X}}^{n} }{\Delta t} &= \frac{ 5 F(\vec{\mathbf{X}}^*,t^{n+1}) + 8 F(\vec{\mathbf{X}}^{n},t^{n}) - F(\vec{\mathbf{X}}^{n-1},t^{n-1})}{12}
\label{colloidEvolution2}
\end{align}

\section{Numerical approach}\label{sec:implementation}
\subsection{Mixed meshless method: Divergence-free velocity and staggered pressure}\label{subsec:discretize}
We present here a new mixed discretization utilizing a divergence free reconstruction for the velocity and a staggered compatible meshless reconstruction for the pressure. The domain is discretized by distributing a set of collocation points $\mathbf{N} := \left\{ x_i \right\}_{i=1,\dots,N_p}$ throughout the domain and along the boundary. We will discuss in Section \ref{subsec:multiresolution} the details of this point distribution, but for now assume that the points are quasi-uniform and provide polynomial unisolvency in the sense of Wendland's text\cite{wendland2004scattered}. For points lying on $\partial \Omega$, we assume an outward facing normal $\hat{\mathbf{n}}_i$ is given. We then seek a set of compactly supported differential operators taking values at each point. To do this, we associate with each point $x_i$ a lengthscale $\epsilon_i$, and construct the $\epsilon_i$-graph of points near $x_i$.
\begin{equation}
\label{epsilonCell}
\mathbf{E}_i = \left\{ (i,j) \text{ such that } ||x_i - x_j|| < \epsilon_i \right\}
\end{equation}
The construction of this neighbor connectivity can be easily performed in $O(N)$ time using standard binning algorithms\cite{wendland2004scattered}. The discretization used here is an extension of the standard MLS approach, in which a differential operator is constructed by finding a local polynomial reconstruction of function values $u_j:=u(x_j)$ as the solution to the weighted $\ell_2$-optimization 
\begin{equation}
\label{standardMLSa}
{q}^* =\argmin\limits_{{q}\in \pi_m}  \left\{ \sum_{j=1}^{N_p} \left[ u_j - q_j   \right]^2 W_{ij} \right\}
\end{equation} 
where $\pi_m$ denotes the $m^{th}$ order polynomials, and $W_{ij}$ denotes a positive weight depending on the relative distances of $x_i$ and $x_j$, satisfying $W_{ij}=W_{ji}$. To approximate an operator $D^\alpha$ evaluated at $x_i$, the operator is applied to the polynomial reconstruction
\begin{equation}
\label{standardMLSb}
D^\alpha u_i \approx D^\alpha_h u_i := D^\alpha q^*(x_i).
\end{equation}  
By choosing $W_{ij}$ with compact support, this efficiently generates finite difference-like stencils that can be applied to complex geometries. A detailed discussion of these problems and their stable solution can be found in the literature, and we refer the reader to the following sources for details\cite{mirzaei2011generalized,wendland2004scattered,kim2003point}. We note that since $\ell_2$-minimization has an analytic solution, Equation \ref{standardMLSb} provides a means of generating finite-difference like stencils of the form
\begin{equation}
\label{standardMLS_FD}
D^\alpha_h u_i = \sum_{j\in supp(W_{ij})} \alpha_{ij} u_j
\end{equation}
where $\alpha_{ij}$ are coefficients obtained from the solution to Equation \ref{standardMLSa}. Details of how these coefficients are calculated can be found in our previous works \cite{stagMLS,compactMLS}. Unlike standard finite difference schemes or Galerkin methods, these coefficients do not inherit the symmetry of the corresponding continuous operators (e.g. when approximating a gradient, $\alpha_{ij} \ne -\alpha_{ji}$ except for special cases of uniform particle arrangements).
\begin{figure}[h!]
  \centering
  \includegraphics[height=2in]{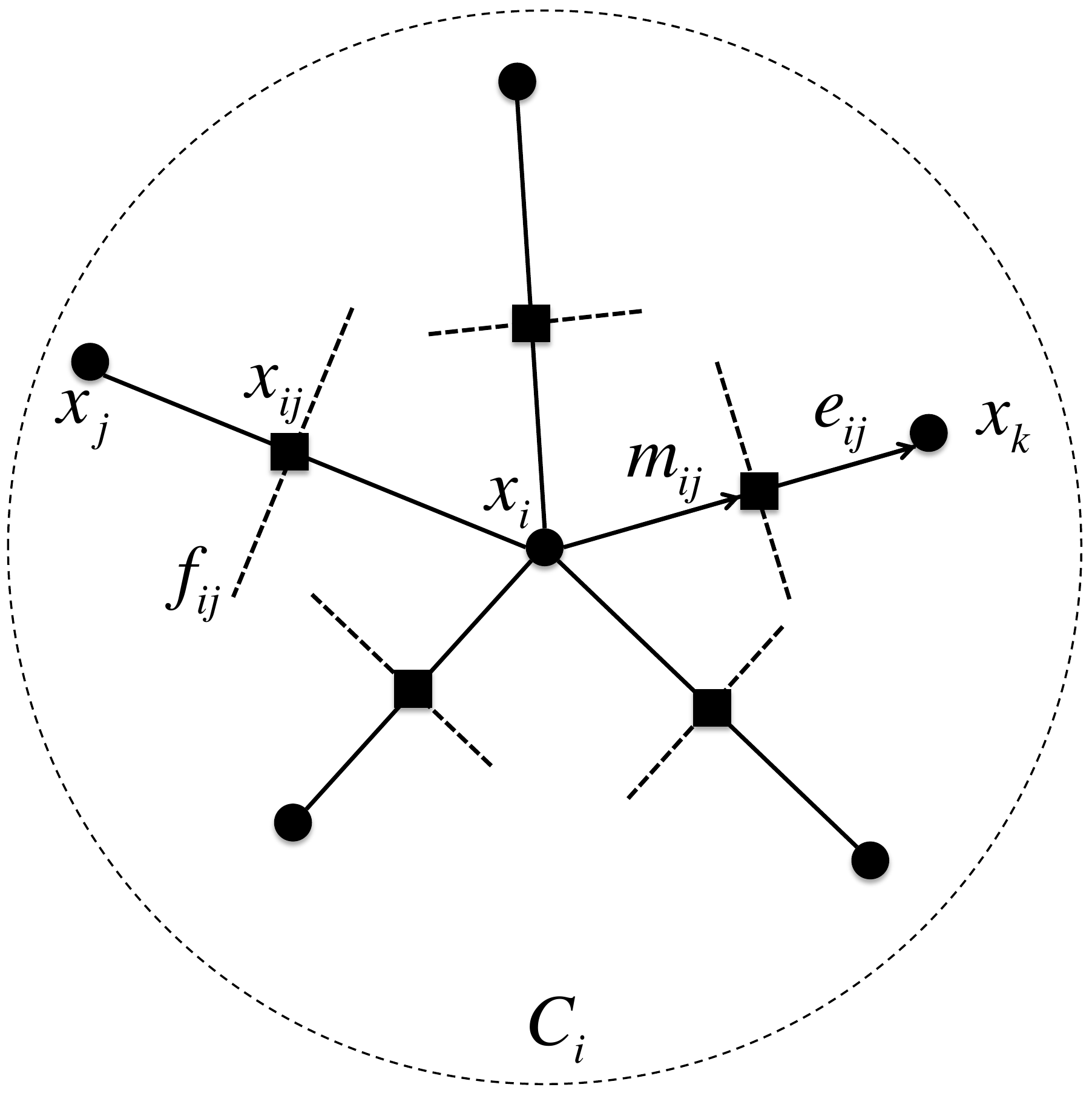} 
  \caption{A local primal-dual grid complex induced by a point $\bm{x}_i$ and its $\varepsilon$-neighborhood $N^\varepsilon_i$. The primal edges $\bm{e}_{ij}$, midpoints $\bm{x}_{ij}$ and mid-edges $\bm{m}_{ij}$ are physical mesh entities. The dual cell $C_i$ and the dual faces $\bm{f}_{ij}$ are virtual mesh entities\cite{stagMLS}.}
  \label{stencil}  
\end{figure}

In the staggered moving least squares method \cite{stagMLS}, for each point $\mathbf{x}_i$ we assemble the $\epsilon_i$-graph $\mathbf{E_i}$ (Figure \ref{stencil}). We identify with each edge a virtual dual face and with each node a virtual cell, and seek a meshless discretization of divergence and gradient operators mimicking the algebraic structure of staggered finite volume methods. Defining $\mathbf{E}:=\cup_i \mathbf{E}_i$, we seek staggered divergence and gradient discretizations of the form
\begin{align}
div_h : \mathbf{E} \rightarrow \mathbf{N}\\
grad_h : \mathbf{N}   \rightarrow \mathbf{E}.
\end{align}

By associating integral degrees of freedom with each edge, a discrete combinatorial gradient is defined
\begin{equation}
grad_h(\phi)_{ij} = \int_{e_{ij}} \nabla \phi \cdot d\mathbf{s} = \phi_j - \phi_i
\end{equation}
that is exact via the fundamental theorem of calculus. 

A divergence operator can then be constructed by first defining the radial component function

\begin{equation}
u_{i\rightarrow}(\bm{x}) =  \bm{u}(\bm{x})\cdot 2\left(\bm{x}-\bm{x}_i\right)
\label{eq:radial}
\end{equation}

When restricted to each edge, this provides the scalar component of the vector field aligned with each edge. This quantity plays a role similar to a flux in the finite volume method, where a vector field is reconstructed from the normal component of fluxes through faces at cell boundaries. 

For a sufficiently smooth function $\mathbf{u}$, the following identities hold relating the radial component function to the vector field and its divergence at the nodes of the graph:
\begin{equation}
\bm{u}(\bm{x}_i) = \frac12\nabla u_{i\rightarrow}(\bm{x}_i)
\quad\mbox{and}\quad
\nabla\cdot\bm{u}(\bm{x}_i) =\frac14 \nabla \cdot \nabla u_{i\rightarrow}(\bm{x}_i) .
\label{eq:radial-lap}
\end{equation}
We can then reconstruct vector fields and their divergences at nodes by sampling their radial component functions at edges. Defining the edge midpoint $\bm{x}_{ij}=\left(\frac{\bm{x}_j+\bm{x}_i}{2}\right)$
\begin{equation}
{q}^* =\argmin\limits_{{q}\in \pi_m}  \left\{ \sum_{j} \left[ \bm{u}(\bm{x}_{ij})\cdot2\left(\bm{x}_{ij}-\bm{x}_i\right) - q\left(\bm{x}_{ij}\right)   \right]^2 W_{ij} \right\}
\end{equation}
discrete approximations to $\bm{u}$ and its divergence at the nodes are then obtained via
\begin{equation}
\bm{u}_h(\bm{x}_i) = \frac12\nabla q^*(\bm{x}_i)
\quad\mbox{and}\quad
\nabla_h\cdot\bm{u}(\bm{x}_i) =\frac14 \nabla \cdot \nabla q^*(\bm{x}_i) .
\end{equation}

If we instead reconstruct from the discrete gradient operator at the edge midpoints
\begin{equation}
{p}^* =\argmin\limits_{{p}\in \pi_m}  \left\{ \sum_{j} \left[ grad_h(\phi)_{ij} - q\left(\bm{x}_{ij}\right)   \right]^2 W_{ij} \right\}
\label{eq:gradlapQP}
\end{equation}
we then obtain the following discrete gradient and divergence operators at each node
\begin{equation}
\nabla_h \phi(\bm{x}_i) = \frac12\nabla p^*(\bm{x}_i)
\quad\mbox{and}\quad
\nabla^2_h\cdot\phi(\bm{x}_i) =\frac14 \nabla \cdot \nabla p^*(\bm{x}_i) .
\label{eq:gradlap}
\end{equation}

For the sake of brevity, we refer to a previous work for details regarding the derivation, analysis and implementation of the SMLS method\cite{stagMLS}. When the operators given by Equation \ref{eq:gradlap} are used to discretize a Poisson equation with Neumann data using an $m^{th}$ order polynomial reconstruction, the resulting solution provides $m^{th}$ order convergence in both a discrete $\ell_2$ and $H^1$ norm. For boundary particles, Neumann boundary conditions are easily handled by adding a constraint to Equation \ref{eq:gradlapQP} to enforce that the reconstructed gradient exactly satisfy the boundary condition at the point $\bm{x}_i$. Similar to Equation \ref{standardMLS_FD}, this process yields a pair of finite difference formulae that can be used to discretize the gradient and Laplacian terms in Equations \ref{StokesEqn2a} and \ref{StokesEqn2b}, respectively.

With the gradient and Laplacian operators discretized, we finally seek a discretization of the viscous operator compatible with the divergence-free constraint. To do this, we introduce a standard MLS discretization defined instead over the space of divergence-free vector polynoimals. we define the set of such $m^{th}$ order polynomialsas $\mathbf{\pi}^{div}_m$. Unlike finite element methods for which the construction of a divergence-free basis is challenging due to continuity constraints, for this setting a basis can easily be constructed, e.g. in the two-dimensional case with $m=1$, for any $\mathbf{\phi}\in \pi^{div}_1$

\begin{equation}
\label{divfreebasis}
\mathbf{\phi} \in span\left\{ \binom{1}{0},\binom{0}{1},\binom{y}{0},\binom{0}{x},\binom{x}{-y}\right\}
\end{equation}
Therefore we can procede with a standard MLS reconstruction, albeit using a modified basis.
\begin{equation}
\label{divfreeMLSa}
\mathbf{p}^* =\argmin\limits_{{p}\in \pi^{div}_m}  \left\{ \sum_{j=1}^{N_p} \left[ \mathbf{u}_j - \mathbf{p}_j   \right]^2 W_{ij} \right\}
\end{equation}
The viscous operator may then be discretized as
\begin{equation}
\label{divfreeMLSb}
\nabla \times \nabla \times \mathbf{u}(x_i) \approx \nabla \times \nabla \times_h \mathbf{u}_i := \nabla \times \nabla \times \mathbf{p}^*(x_i)
\end{equation}
If we denote the $\mathbb{R}^{dim(\mathbf{\pi}^{div}_m)\times d}$ matrix containing the polynomial basis functions in Equation \ref{divfreebasis} as $\mathbf{P}(x)$, Equation \ref{divfreeMLSb} can explicitly be written as
\begin{equation}
\nabla \times \nabla \times_h \mathbf{u}(x_i) = \nabla \times \nabla \times \mathbf{P}(x_i)^\intercal  \mathbf{M}(\mathbf{x}_i)^\dagger \sum_{j=1}^{N_p} \mathbf{P} (x_j) W_{ij}  \mathbf{u}_j ,
\end{equation}
where $\mathbf{M}^\dagger$ denotes the pseudo-inverse of the weighted normal equations matrix
\begin{equation}
\mathbf{M}(\mathbf{x}_i) = \sum_{j=1}^{N_p} \mathbf{P} (x_j) W_{ij} \mathbf{P}^\intercal (x_j)
\end{equation}
and again we obtain a finite difference-like formula similar to Equaton \ref{standardMLS_FD} that may be used to discretize the viscous operator.

As mentioned in Section \ref{sec:Introduction}, the lack of a mesh or symmetry in particle arrangement makes it difficult to pose a rigorous analysis as to the sense in which this reconstruction is actually faithful to the divergence-free constraint. It is unclear whether a discrete mass conservation statement can be proven here. Unlike Wendland's approach using matrix-valued RBFs to obtain an analytically divergence-free velocity space \cite{wendland2009divergence}, our approach is only divergence-free in the local polynomial reconstruction. In Huerta's work \cite{huerta2004pseudo}, a similar reconstruction was used but for a weak formulation, allowing them to demonstrate that the reconstruction passes a numerical inf-sup test. We will demonstrate in Section \ref{subsec:wannier} that the use of this divergence-free reconstruction similarly provides stable approximation for this strong form discretization while maintaining high-order convergence for both pressure and velocity, and leave a rigorous analysis for a future work.

To summarize, similar to Equation \ref{standardMLS_FD}, Equations \ref{eq:gradlap} and \ref{divfreeMLSb} provide the following set of finite-difference like stencils 
\begin{align}
\label{MMMstencils}
\nabla^2_h p_i = \sum_{j \in supp(W_{ij})} \alpha^1_{ij} p_j\\
\nabla_h p_i = \sum_{j \in supp(W_{ij})} \alpha^2_{ij} p_j\\
\nabla \times \nabla \times_h \mathbf{u}_i = \sum_{j \in supp(W_{ij})} \alpha^3_{ij} \mathbf{u}_j\\
\end{align}
where $\alpha^1,\alpha^2,\alpha^3$ are coefficients obtained from the solution of each optimization problem. These stencils can then be used to discretize the Stokes equations in a manner identical to standard finite-difference methods: a global matrix is assembled by discretizing Equations \ref{StokesEqn2a} and \ref{StokesEqn2b} using Equation \ref{MMMstencils} at each point in the interior of the domain. For points lying on the boundary, the Dirichlet boundary condition in Equation \ref{StokesEqn2a} is enforced by adding a one to the diagonal of the global matrix and setting the imposed value on the right-hand side. The Neumann boundary conditions for the pressure in Equation \ref{StokesEqn2b} are enforced implicitly via a constraint in Equation \ref{eq:gradlapQP}, following the procedure given by \citet{stagMLS}. As in any method discretizing the Poisson problem with Neumann boundary conditions, the discrete Poisson operator contains a constant vector in its right null space that must be accounted for in order to obtain a unique solution\cite{strikwerda1984finite,bochev2005finite}. In this work we handle this by adding a single Lagrange multiplier to the Poisson problem to enforce a zero mean pressure.
\begin{equation}
\label{nullSpace}
\sum_{i=1,\dots,N_p} p_i = 0
\end{equation}

\subsection{Multiresolution approach: kernel adaptivity and particle distribution}\label{subsec:multiresolution}
To simulate the suspension flow described in Equation \ref{StokesEqn1} from first principles, it is necessary to resolve both the singular pressure occuring in the narrow lubrication gap and the larger length scales of the particle/boundary interactions that drive the flow. To do this we adaptively select the optimization weights $W_{ij}$ to match the surrounding particle refinement. In standard applications of MLS, the weights are chosen by selecting a positive symmetric kernel with support $\epsilon$. 

\begin{equation}
W_{ij} = W_\epsilon(||x_i - x_j||)
\end{equation}

In this work, we consider kernels of the form $W_\epsilon(r) = (1 - \frac{|r|}{\epsilon})^p$, with $p=4$. For a distribution of particles characterized by a quasi-uniform lengthscale $\Delta x$, the support can be scaled as $\epsilon = C \cdot \Delta x$, where $C$ is selected to ensure that enough neighbors lie in the support of $W_{ij}$ so that Equation \ref{standardMLSa} has a unique solution. Assuming that the points maintain polynomial unisolvency, this can loosely be considered as requiring that the number of neighbors be larger than the dimension of the polynomial space
\begin{equation}
\label{adaptiveReq}
\# \left\{ j \in supp(W_{ij}) \right\} > dim(\pi_m)
\end{equation}
For example, taking $C = m+1$, where $m$ is the order of the MLS reconstruction, provides a slightly conservative interaction radius that is appropriate near corners and boundaries of the domain where one-sided stencils are generated.

To obtain a multiresolution discretization, the support of the kernel is chosen adaptively as a function of position, i.e.  $\epsilon_i := \epsilon(x_i)$. An initial lengthscale characterizing the maximal discretization lengthscale is calculated
\begin{equation}
\Delta x^{\infty} = sup_i \min_{j \ne i} ||x_i - x_j||
\end{equation}
and a preliminary support radius is calculated as $\epsilon_i^{prelim} = (m+1)\Delta x^\infty$. Neighbor lists are then assembled and sorted in order of increasing distance. In light of Equation \ref{adaptiveReq}, the distance to neighbor number $dim(\pi_m)$ is calculated. This distance $r_{ij}^{min}$ denotes the minimum distance to obtain enough neighbors to reconstruct a polynomial under the assumption that all neighbors are $\pi_m$-unisolvent. To account for the possibility of degenerate particle arrangements, the particle support is conservatively chosen as $\epsilon_i = 1.5 r_{ij}^{min}$. The weights used in $\ell_2$ optimization are then defined as
\begin{equation}
W_{ij} = W_{\epsilon_i}(||x_i - x_j||) + W_{\epsilon_j}(||x_i - x_j||)
\end{equation}
These weights are again positive and symmetric in $i$ and $j$, but scale appropriately between regions of high and low resolution. The neighbor lists are then recalculated to match the support of $W_{ij}$, and the global matrix can be assembled. 

Before weights can be calculated however, particles must be distributed throughout the domain. For simplicity, we pursue the same refinement strategy for all cases in this work. The discretization is characterized by a lengthscale $\Delta x^\infty$, the number of refinement levels $M_r$ and the number of refinement layers $M_l$. Each refinement level $L_i$ is characterized by a lengthscale $\Delta x_i = \Delta x^\infty 2^{i-M_r}$. Particles are first placed along colloid boundaries with separation $\Delta x_0$. Then for $i=1,\dots,M_r$, $M_l$ layers of thickness and arclength $\Delta x_i$ are placed marching outward from the boundary. If any of these particles intersect a more highly refined layer belonging to another particle, they are not added. After these nested layers are generated, a Cartesian grid of particles with spacing $\Delta x^{\infty}$ is placed covering the remainder of the domain, and finally the outter boundary is discretized by placing particles around the perimeter with spacing $\Delta x^{\infty}$. For the results presented later, each discretization will be described by the triplet $\{\Delta x^\infty,M_r,M_l\}$. While this refinement strategy is adopted here to allow a simple, consistent description of the refinement in each benchmark, more sophisticated refinement strategies can be easily adapted, provided they ensure sufficient neighbors to perform the MLS reconstruction.

\subsection{Evaluation of viscous stresses at colloid boundaries}\label{subsec:forceCalc}

In order to couple the motion of the colloids to the flow solver via Equation \ref{forcetorque}, it is necessary to derive an approximation to the integral of the stress tensor in Equation \ref{forcetorque}. We describe the process for the force balance here, but the process for the torque balance is identical. To this end, we identify each boundary particle $x_j \in \partial \Omega_i$ with a portion of arclength $d\mathbf{\theta}_j$ and a unit vector pointing outward from the colloid $\vec{\mathbf{n}}_j$. The integral in Equation \ref{forcetorque} can then be partitioned
\begin{equation}
\label{integralChopped}
  \int_{\partial \Omega_i} \mathbf{\sigma}_V \cdot d\mathbf{A} = \sum_{j \in \partial \Omega_i} \int_{d\mathbf{\theta}_j} \mathbf{\sigma}_V \cdot d\mathbf{A}
\end{equation}
A high-order quadrature rule can then be derived by posing the quadrature via a standard MLS reconstruction
\begin{equation}
  \int_{d\mathbf{\theta}_j} \mathbf{\sigma}_V \cdot d\mathbf{A} \approx   \int_{d\mathbf{\theta}_j} \mathbf{\sigma}^* \cdot d\mathbf{A}
\end{equation}
where $\mathbf{\sigma}^*$ is calculated by taking the reconstructions of pressure and velocity
\begin{align}
q^* =\argmin\limits_{q\in \pi_m}  \left\{ \sum_{j=1}^{N_p} \left[ p_j - q_j   \right]^2 W_{ij} \right\}\\
\mathbf{v}^* =\argmin\limits_{\mathbf{v}\in \pi^{div}_m}  \left\{ \sum_{j=1}^{N_p} \left[ \mathbf{u}_j - \mathbf{v}_j   \right]^2 W_{ij} \right\}\\
\end{align}
and calculating
\begin{equation}
\mathbf{\sigma}^*_V = -q^* \mathbf{I} + \frac{\nu}{2}\left( \nabla \mathbf{v}^* + \nabla {\mathbf{v}^*}^\intercal\right)
\end{equation}
which can then be used to evaluate the integral in Equation \ref{integralChopped}. Although accurate to high order, for the cases considered in this work, the use of the following midpoint quadrature rule
\begin{equation}
  \mathbf{F}_d \approx \sum_{j \in \partial \Omega_i} \mathbf{\sigma}^*(x_j) \cdot \vec{\mathbf{n}}_j R_i d\mathbf{\theta} 
\end{equation}
provided no discernable difference, and is used for the results presented. 

\subsection{Efficient solution of the monolithic system}\label{subsec:solvers}

If the colloid motion is coupled together with the flow equations through Equations \ref{StokesEqn2a},\ref{StokesEqn2b} and \ref{forcetorque} and the discretization described in Sections \ref{subsec:discretize} and \ref{subsec:forceCalc} is applied, a monolithic block structured matrix is obtained coupling together velocity, pressure and the degrees of freedome describing the rigid body rotation of each particle. For the purposes of preconditioning, we partition the matrix into the following $2 \times 2$ block matrix system
\begin{equation}\label{eqn:2by2-mat}
\begin{bmatrix}
\mathbf{K} & \mathbf{G} \\
\mathbf{B} & \mathbf{L}
\end{bmatrix}
\begin{bmatrix}
\mathbf{\tilde{u}} \\
p
\end{bmatrix}
 = 
 \begin{bmatrix}
 \mathbf{f} \\
 g
 \end{bmatrix},
\end{equation}
where we combine the variables pertaining to the velocity and particle degrees of freedom into a vector $\bm{\tilde{u}}$, and the blocks are given as follows: $\mathbf{K}$ contains contributions from the viscous operator and viscous coupling through the stress balance at the colloid surface, $\mathbf{G}$ corresponds to the contribution of pressure to the momentum equation and pressure drag in the colloid stress balance, $\mathbf{L}$ corresponds to the Poisson operator, and $\mathbf{B}$ is non-zero only at the boundary and contains the contribution from the $\nu \hat{n} \cdot \nabla \times \nabla \times$ boundary condition. While a formal analysis of these operators is not possible, numerical computation of eigenvalues for small systems suggest that both $\mathbf{K}$ and $\mathbf{L}$ are positive-definite. 

In order to solve the coupled linear system, the preconditioned general minimal residual (PGMRes) method is used.  As is well-known, an efficient and robust preconditioner is essential for the overall performance of the PGMRes method.  In this work, we adopt the triangular block preconditioner based on the following block factortization of the $2 \times 2$ system \eqref{eqn:2by2-mat}.   
\begin{equation}\label{eqn:block-fac}
\begin{bmatrix}
\mathbf{K} & \mathbf{G} \\
\mathbf{B} & \mathbf{L}
\end{bmatrix}
= 
\begin{bmatrix}
\mathbf{S} & \mathbf{G} \\
\mathbf{0} & \mathbf{L}
\end{bmatrix}
\begin{bmatrix}
\mathbf{I} & \mathbf{0}\\
\mathbf{L}^{-1} \mathbf{B} & \mathbf{I}
\end{bmatrix},
\end{equation}
where $\mathbf{S} = \mathbf{K} - \mathbf{G} \mathbf{L}^{-1} \mathbf{B}$ is the Schur complement.  As usual, explicitly assembling the Schur complement $\mathbf{S}$, however, involves inverting $\mathbf{L}$ explicitly, which is expensive and the resulting $\mathbf{S}$ is dense.  Therefore, we approximate $\mathbf{L}$ by its diagonal entries, i.e., $\text{diag}(\mathbf{L})$.  This gives an approximate Schur complement $\widetilde{\mathbf{S}} = \mathbf{K} - \mathbf{G} \; \text{diag}(\mathbf{L})^{-1} \mathbf{B}$ and the corresponding block triangular preconditioner is given by
\begin{equation}\label{eqn:prec}
\begin{bmatrix}
\widetilde{\mathbf{S}}  & \mathbf{G} \\
\mathbf{0} & \mathbf{L}
\end{bmatrix}.
\end{equation}  
We want to comment that such a block preconditioner is different from the traditional block preconditoners designed for Stokes problems that involve the other Schur complment $\mathbf{L} - \mathbf{B} \mathbf{K}^{-1} \mathbf{G}$ and its approximation $\mathbf{L} - \mathbf{B} \; \text{diag}(\mathbf{K})^{-1} \mathbf{G}$.  The reason is that, in our setting, $\mathbf{L}$ corresponds to the discretization of $\nabla^2$. Numerical experiments suggest that it is close to being diagonally dominated and, therefore, can be well-approximated by its diagonal.  However, $\mathbf{K}$ comes from $\nu \nabla \times \nabla \times$ and numerical experiments showed that it was far from diagonally dominant.  Therefore, the standard choice of $\text{diag}(\mathbf{K})$ in approximating the Schur complement of the upper left block proved not to provide a good approximation. 

In order to apply the block preconditioner \eqref{eqn:prec}, both $\mathbf{L}$ and $\widetilde{\mathbf{S}}$ must be inverted.  In our implementation, we use the AMG method to efficiently solve them both. Numerical benchmarking of the block preconditioner is provided in the following section. 

\section{Stokes solver: Validation, performance and convergence study}\label{sec:stokesval}
To validate the mixed scheme, we begin by studying cases where a single colloid has a known position, velocity, and angular velocity. This will allow a study of the convergence of the Stokes solver and force/torque calculation before coupling colloid dynamics.
\subsection{Wannier flow}\label{subsec:wannier}
Wannier flow consists of two cylinders of radii $R_1$ and $R_2$ with eccentricity $e$ rotating with fixed angular velocity $\omega_1$ and $\omega_2$. Originally studied as a model for lubrication forces in journal bearings\cite{wannier1950contribution}, the flow's analytic solution for velocity, pressure, and the force on the inner cylinder have made this a standard bechmark in validating high-order methods in the presence of curvilinear geometry\cite{karniadakis2013spectral}.

\begin{figure}[h!]
  \centering

\includegraphics[width=1.0\textwidth]{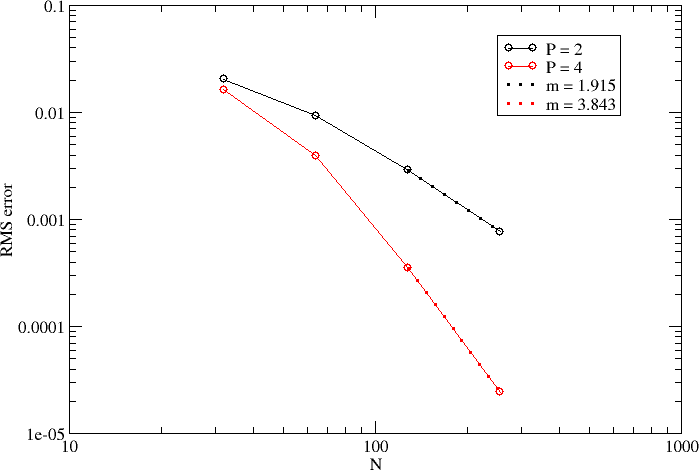}
\includegraphics[width=1.0\textwidth]{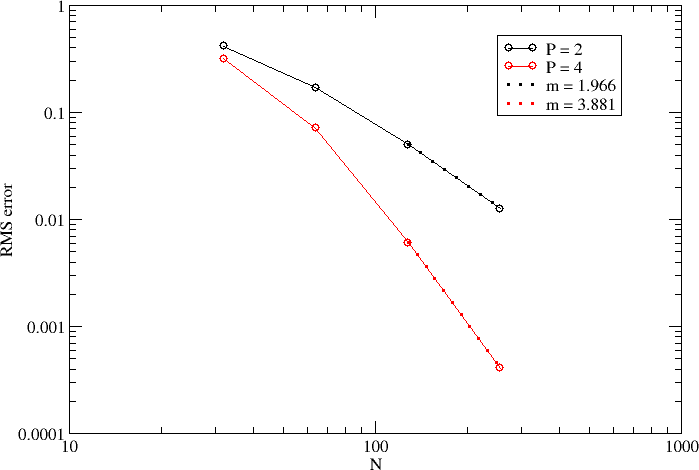}

  \caption{RMS error for velocity (top) and pressure (bottom). Second-order discretizations are given in black, while fourth-order are given in red. Regession of the form $y=cN^{-m}$ demonstrates equal-order optimal convergence for both velocity and pressure. }
  \label{wannier1}
\end{figure}

For this case, we first consider a family of discretizations characterized by the particle distribution $\{N^{-1},1,1\}$, meaning that no adaptive refinement is applied. The parameter $N$ represents the number of particles distributed in each direction, and $N^{-1}$ provides the grid spacing between particles. We begin by considering cylinders of radii $R_1 = \pi/2$, $R_2 = \pi/10$, with eccentricity $\pi/5$, and impose as angular velocities $\dot{\mathbf{\Theta}}_1 = 1/\pi$ and $\dot{\mathbf{\Theta}}_2 = 10/\pi$. We demonstrate that, for an $m^{th}$-order reconstruction of the velocity and pressure, equal order optimal $m^{th}$-order convergence is obtained for both the velocity and the pressure in Figure \ref{wannier1}. The results for all of these simulations take less than a minute to run on a desktop computer with no hardware acceleration.

In Table \ref{tbl:block-prec-stokes} we present benchmarks of the preconditioner for this problem. We that the number of iterations roughly remain the same when the total number of degrees of freedom increases which leads to nearly linearly growth in total CPU time, demonstrating the near optimality and robustness of the block preconditioner.

\begin{table}[htp]
\begin{center}
\begin{tabular}{|c|c|c|c|c|}
\hline 
Degrees of  & Number of  & Setup Time  & Solve Time  & Total Time  \\
freedom & Iterations & (second) & (second) & (second) \\
\hline
 2,715   & 4    & 0.169  & 0.254  & 0.423 \\ 
10,047  & 7   & 0.378  & 1.014  & 1.392 \\
22,026    & 8   & 0.770  & 2.450  & 3.220 \\
38,665    & 8   & 1.221  & 5.758   & 6.979 \\
59,943    & 7   & 1.715  & 6.062   & 7.777 \\
\hline
\end{tabular}
\end{center}
\caption{Preconditioner performance for Stokes problem (PGMRes with stopping criterion that relative residual is less than $10^{-6}$)}
\label{tbl:block-prec-stokes}
\end{table}%
 
\begin{figure}[h!]
  \centering
  \includegraphics[width=0.7\textwidth]{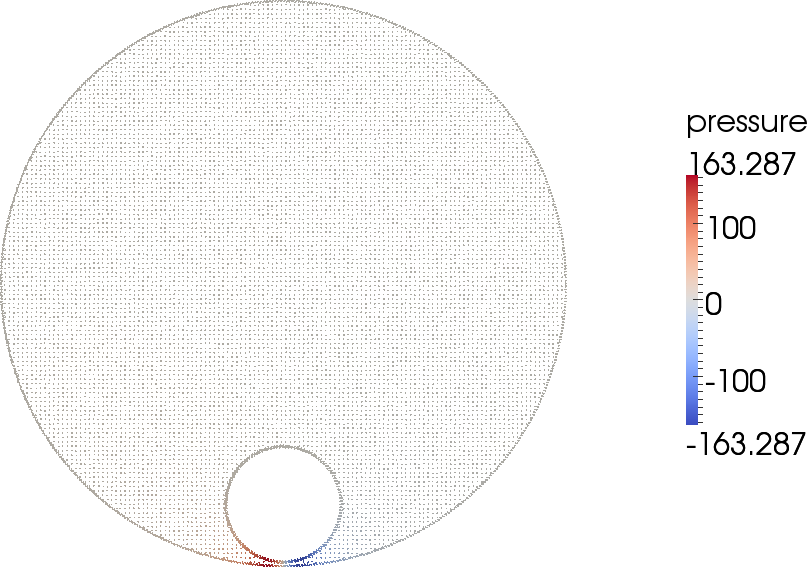}
  \caption{Nearly singular pressure distribution for gap width $a_1/32$.}
  \label{wannier2}
\end{figure}

\begin{figure}[h!]
  \centering
  \includegraphics[width=1.0\textwidth]{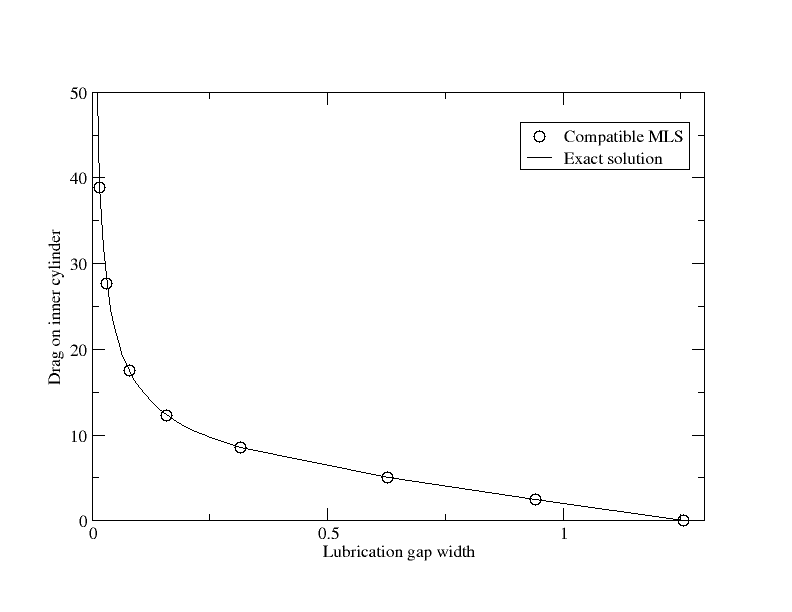}
  \caption{Accurate resolution of singular lubrication force as gap width is reduced.}
  \label{wannier3}
\end{figure}

Next we systematically study the singular lubrication limit as the cylinder separation approaches zero. To do this, we use a fourth-order velocity and pressure reconstruction and apply a mild refinement at around each cylinder using the family of discretizations $\{N^{-1},2,2\}$ and choose $N$ sufficiently large so that there are always at least four points spanning the gap separating the two cylinders. This allows just enough resolution that the stencils within the gap maintain polynomial unisolvency.  The resulting geometry and pressure distribution is shown in Figure \ref{wannier2}, and we demonstrate excellent agreement with the theoretical solution in Figure \ref{wannier3}. While more aggresive adaptivity is certainly possible, we adopt this strategy to demonstrate heuristically that roughly five particles per smallest lengthscale is adequate to obtain accurate results. For the remainder of this work, we use this heuristic to select discretization parameters.

\subsection{Flow past cylinder in channel}

\begin{figure}[h!]
  \centering
  \includegraphics[width=1.0\textwidth]{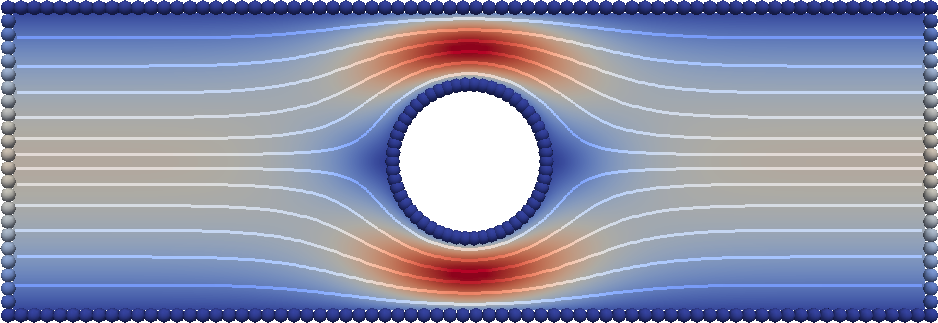}
  \includegraphics[width=1.0\textwidth]{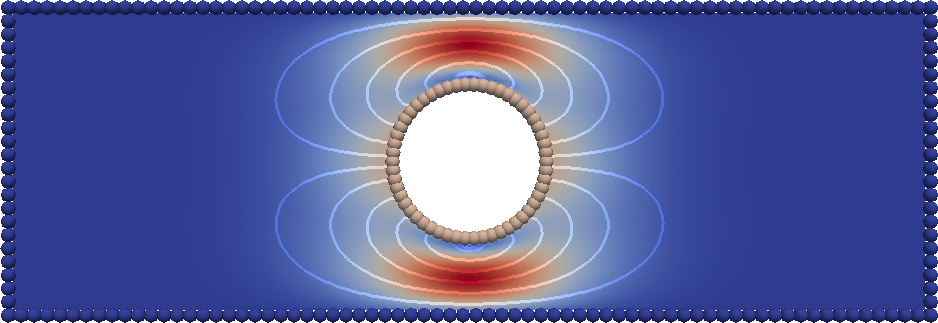}
  \caption{Velocity magnitude and streamlines for the cases: $V=0$, $U=1$, $a=0.5$ (top) and $V=1$, $U=0$, $a=0.5$ (bottom).}
  \label{JY0a}
\end{figure}

\begin{figure}[h!]
  \centering
  \includegraphics[width=1.0\textwidth]{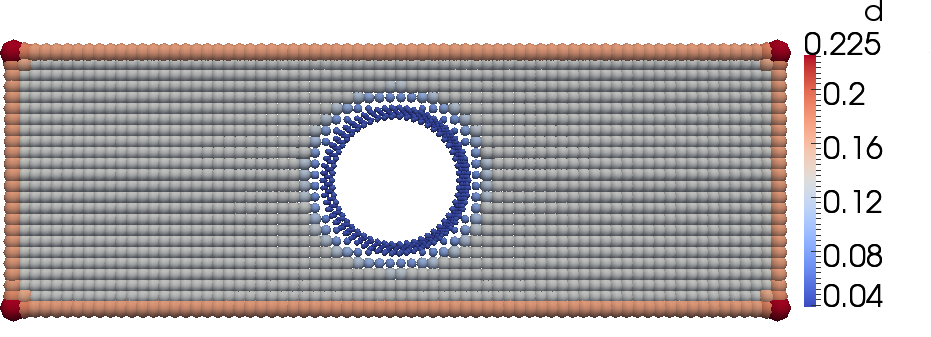}
  \caption{Particle adaptivity visualized by rendering spheres at each point with diameter proportional to $\epsilon_i$.}
  \label{JY0b}
\end{figure}

\begin{figure}[h!]
  \centering
  \includegraphics[width=1.0\textwidth]{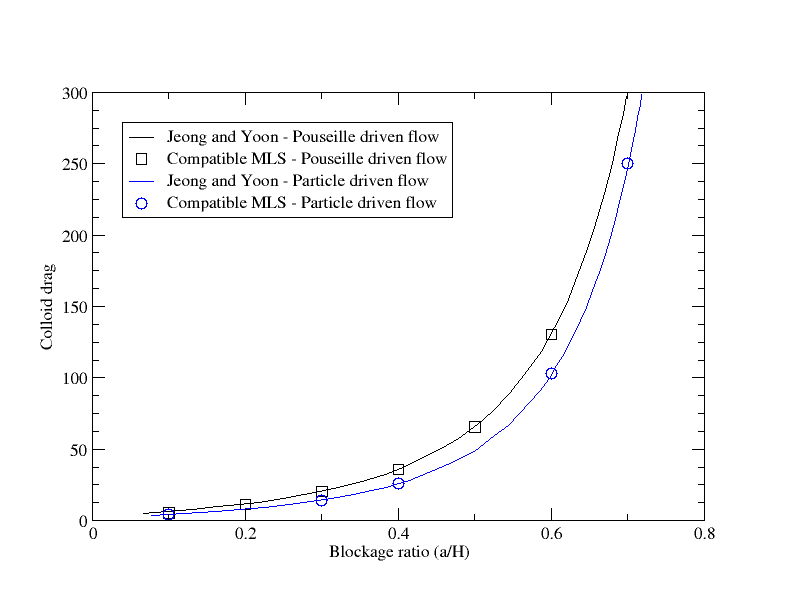}
  \caption{Force for flow driven by Poiseuille flow and by colloid with comparison to exact solution\cite{jeong2014two}}
  \label{JY1}
\end{figure}

\begin{figure}[h!]
  \centering
  \includegraphics[width=1.0\textwidth]{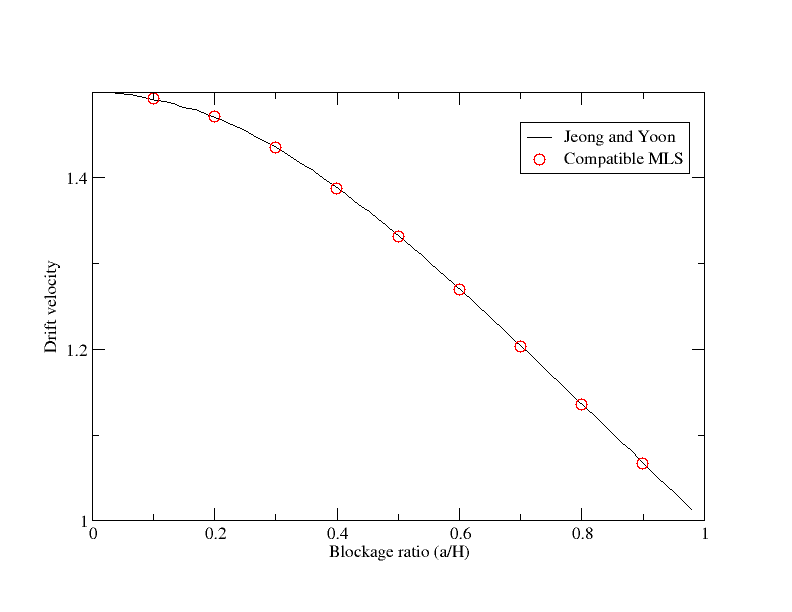}
  \caption{Drift velocity when zero force/torque condition is used to couple colloid motion to flow with comparison to exact solution\cite{jeong2014two}}
  \label{JY2}
\end{figure}

To validate the accuracy of the discretization when the colloid dynamics are coupled to the Stokes solver, we consider flow past a single cylinder of radius $a$ in a channel of length $L$ and height $H$ (See Figures \ref{JY0a},\ref{JY0b}). We consider two cases: one in which zero flow is imposed at the boundary and the particle moves with a constant velocity $V$, and another in which the particle is fixed and a Pouseille flow of magnitude $U$ is imposed at the inlet and outlet. In the limit $L \rightarrow \infty$ an exact solution exists for the drag exerted on the particle as a function of blockage ratio $a/H$. By linearity, \citet{jeong2014two} use this relation to deduce the drift velocity of a particle free to advect under the flow by balancing the flow induced drag with the drag induced by the particle motion. We will calculate this by solving the monolithic system of Equations \ref{StokesEqn2a},\ref{StokesEqn2b} and \ref{forcetorque}. For this and the remainder of this work we use a fourth-order reconstruction. We adopt a $\{N^{-1},3,1\}$ family of refinement and, for a given blockage ratio, refine until approximately five points span the gap between the particle and the wall. We consider the choice of parameters $H=2$, $L=6$, and vary the particle radius $a$. We demonstrate in Figure \ref{JY1} that the corresponding drag forces are very well reproduced, even as the blockage ratio approaches $1$ and lubrication effects become dominant. Finally, having validated that drag forces are accurately resolved, we couple the colloid dynamics together with the Stokes solver via Equation \ref{forcetorque} to match Jeong and Yoon's relation for drift velocity (Figure \ref{JY2}).


\section{Results: large boundary deformation}\label{sec:results}
\begin{figure}[h!]
  \centering
  \includegraphics[width=0.5\textwidth]{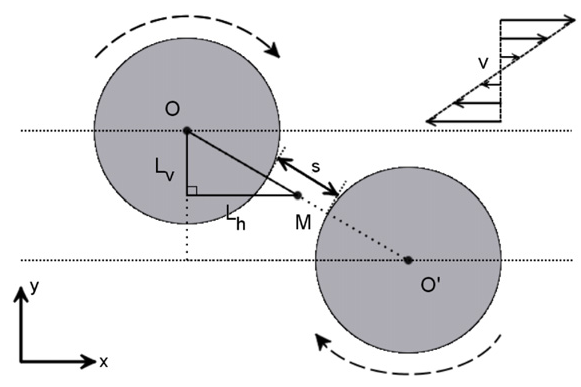}
  \caption{Geometry description for shear flow case.}
  \label{shearFlow}
\end{figure}

\begin{figure}[h!]
  \centering
  \includegraphics[width=1.0\textwidth]{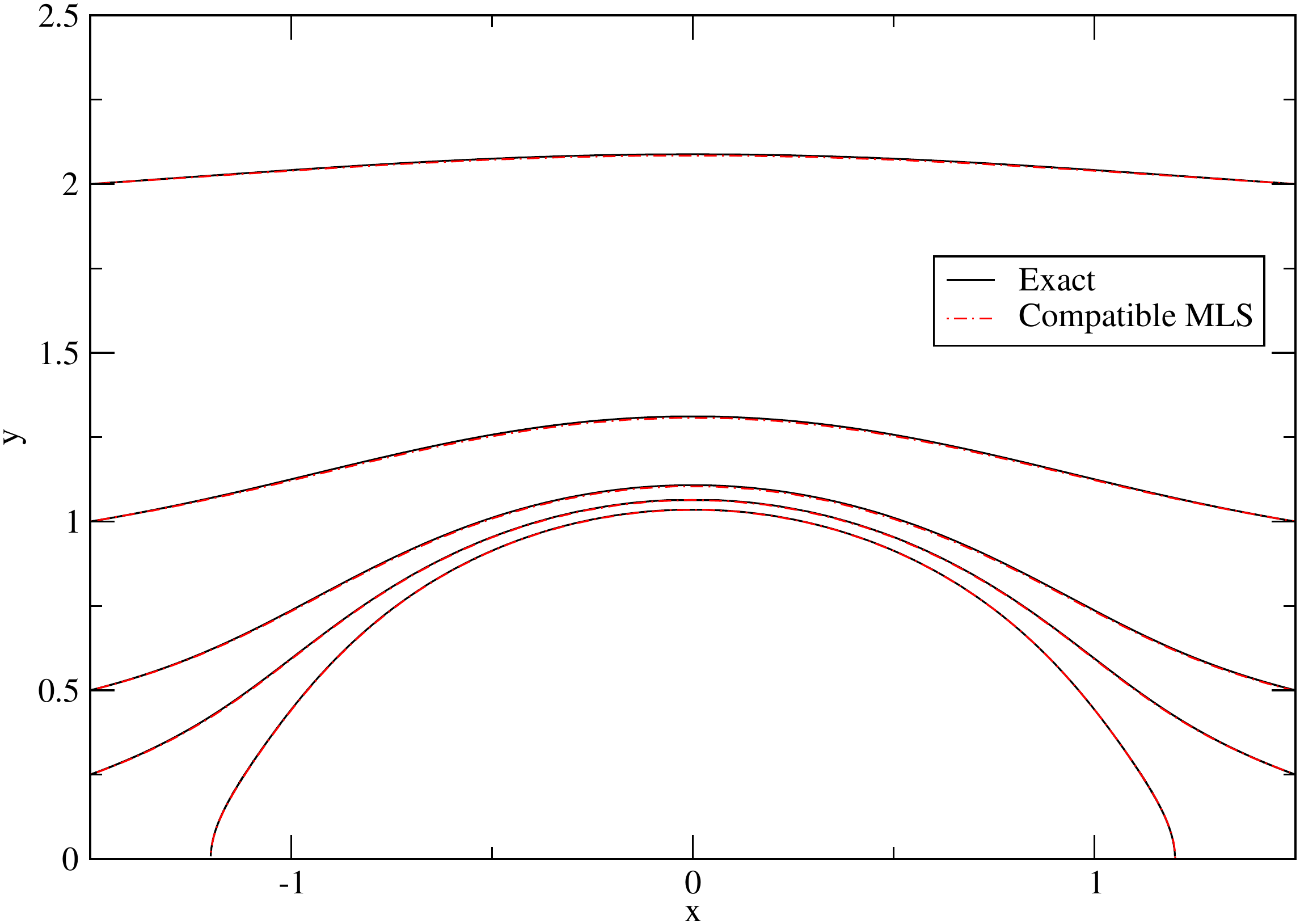}
  \caption{Trajectory of top left particle for varying choice of initial colloid configurarion. For closed orbit, separation is driven by lubrication forces occuring in gap width of size $a/50$.}
  \label{trajectory}
\end{figure}

\begin{figure}
\centering
\begin{subfigure}{.5\textwidth}
  \centering
  \includegraphics[width=.9\linewidth]{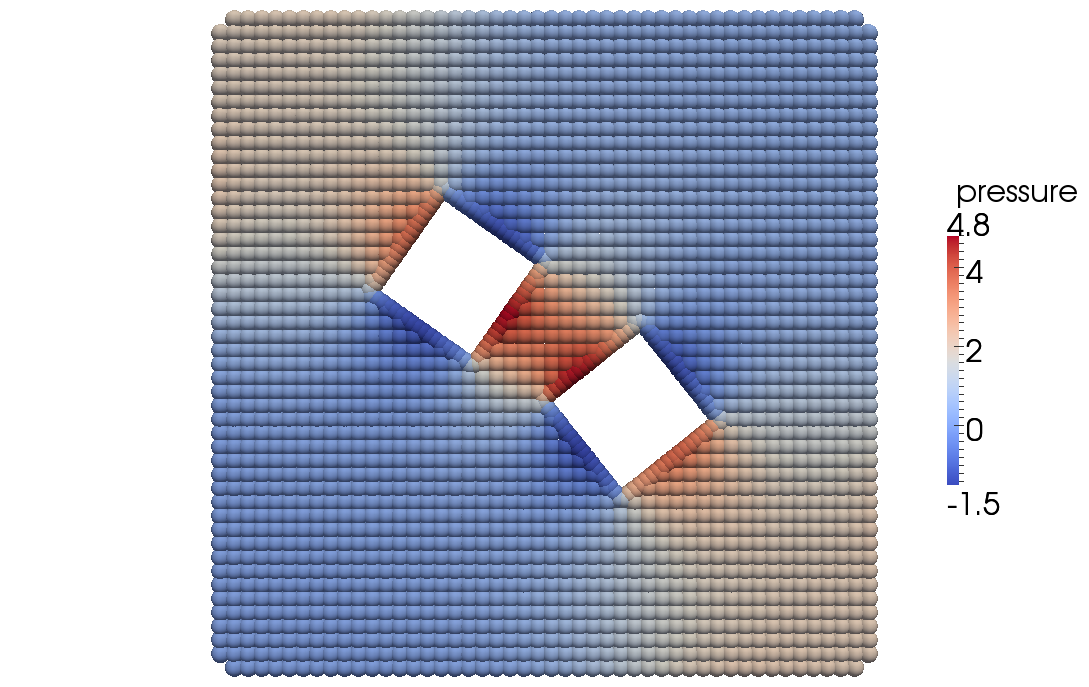}
  \label{fig:sub1}
\end{subfigure}%
\begin{subfigure}{.5\textwidth}
  \centering
  \includegraphics[width=1.10\linewidth]{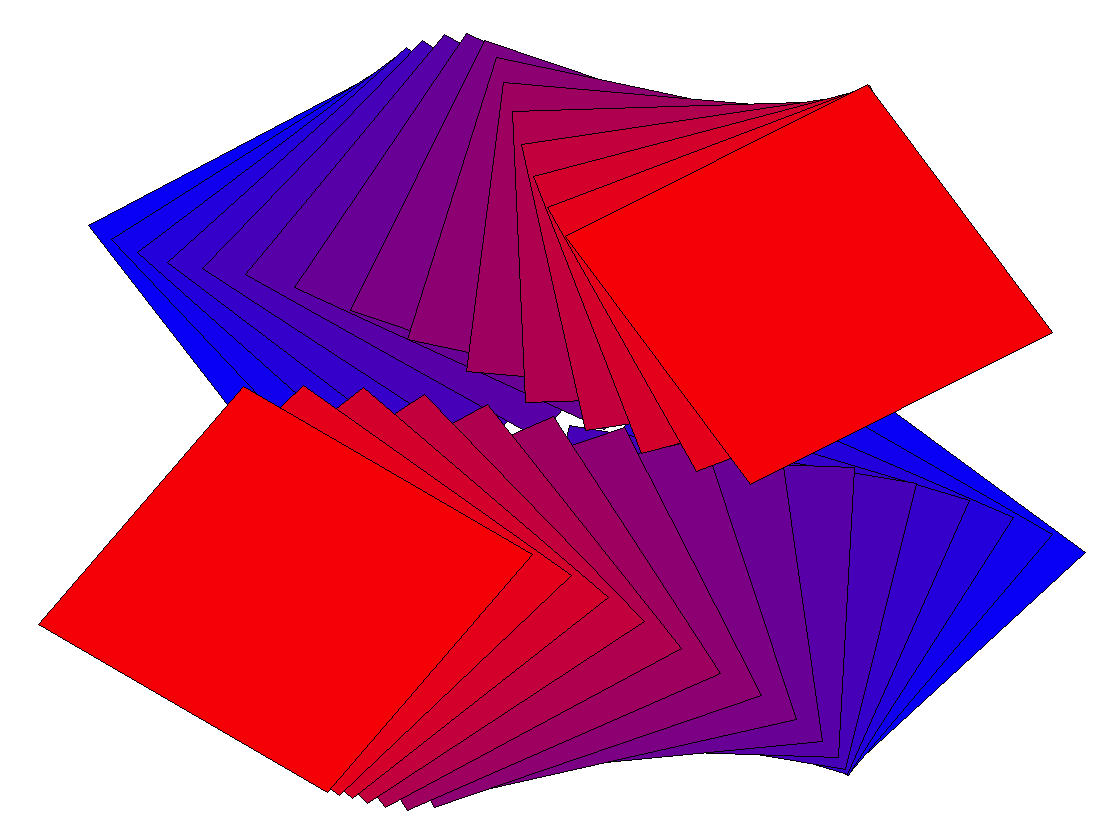}
  \label{fig:sub2}
\end{subfigure}
\caption{(Left) Initial configuration and pressure distribution for square particles of unit width in shear flow. (Right) Trajectory of square particles colored by initial position (red) to final position (blue).}
\label{fig:test}
\end{figure}

\begin{figure}
  \centering
  \includegraphics[width=0.6\linewidth]{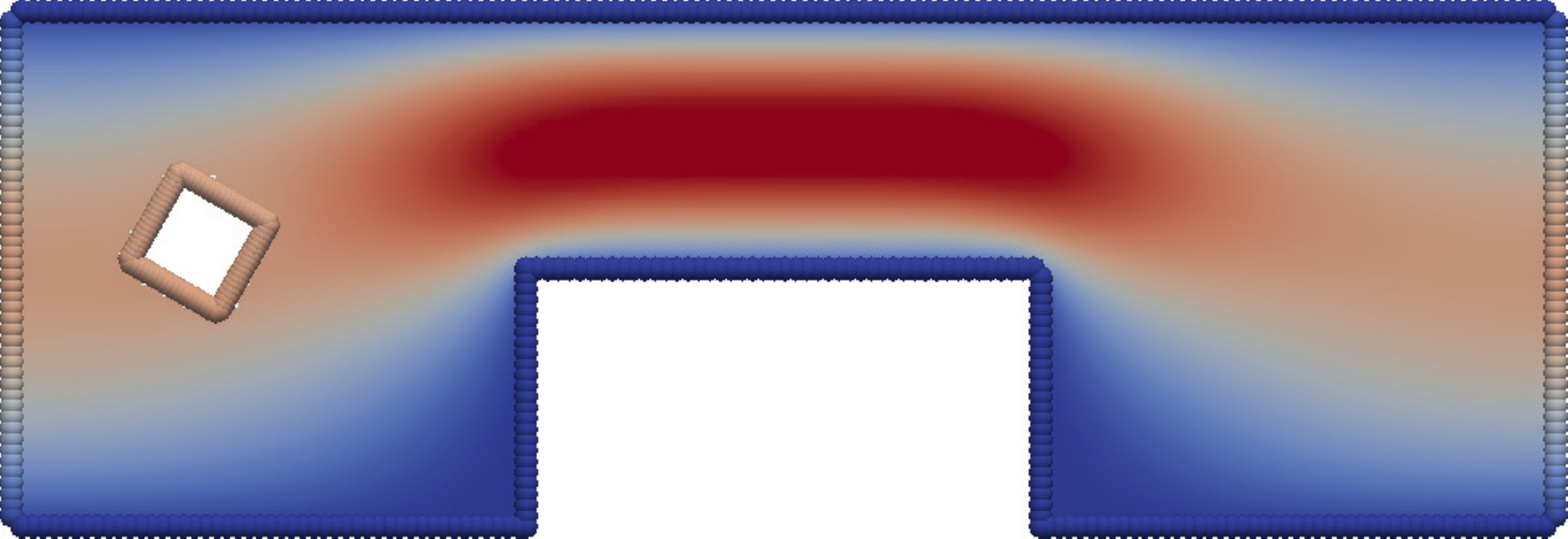}
  \includegraphics[width=0.6\linewidth]{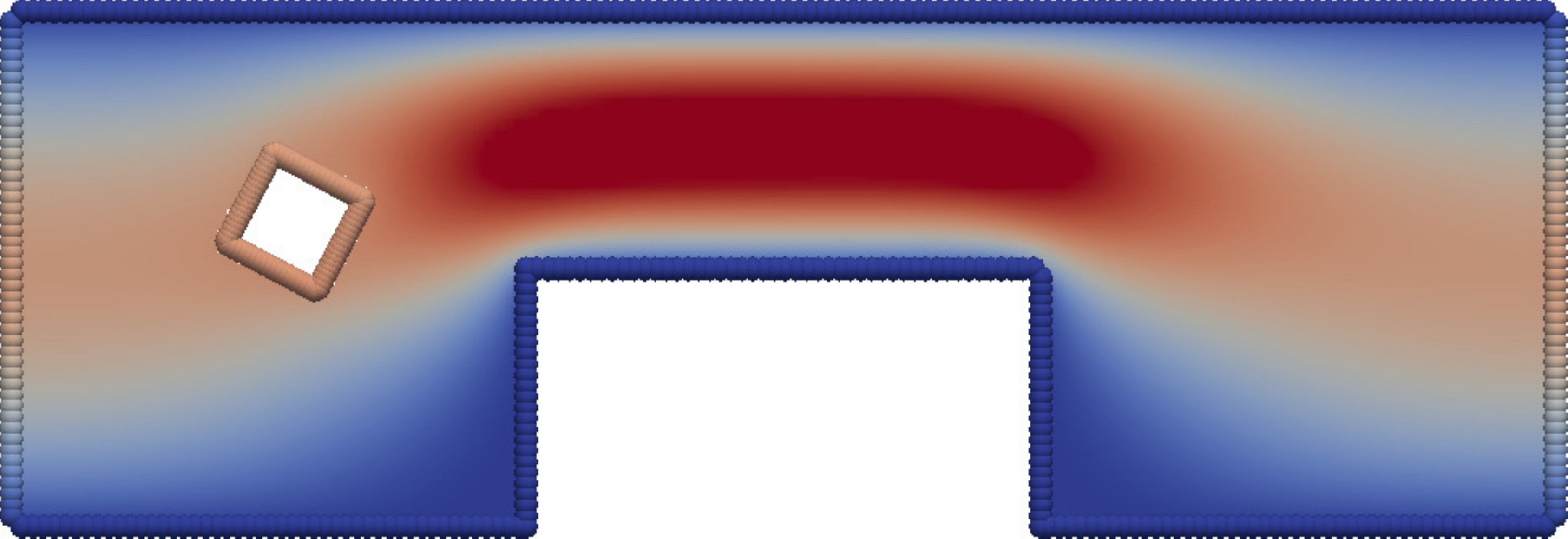}
  \includegraphics[width=0.6\linewidth]{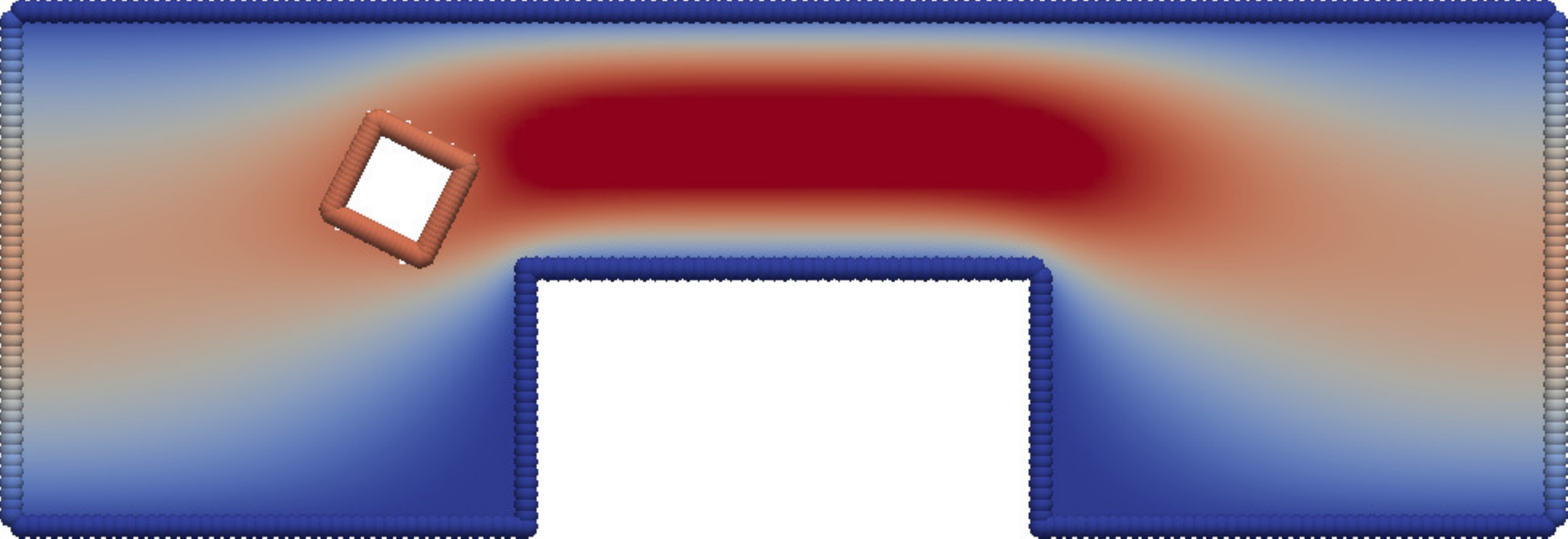}
  \includegraphics[width=0.6\linewidth]{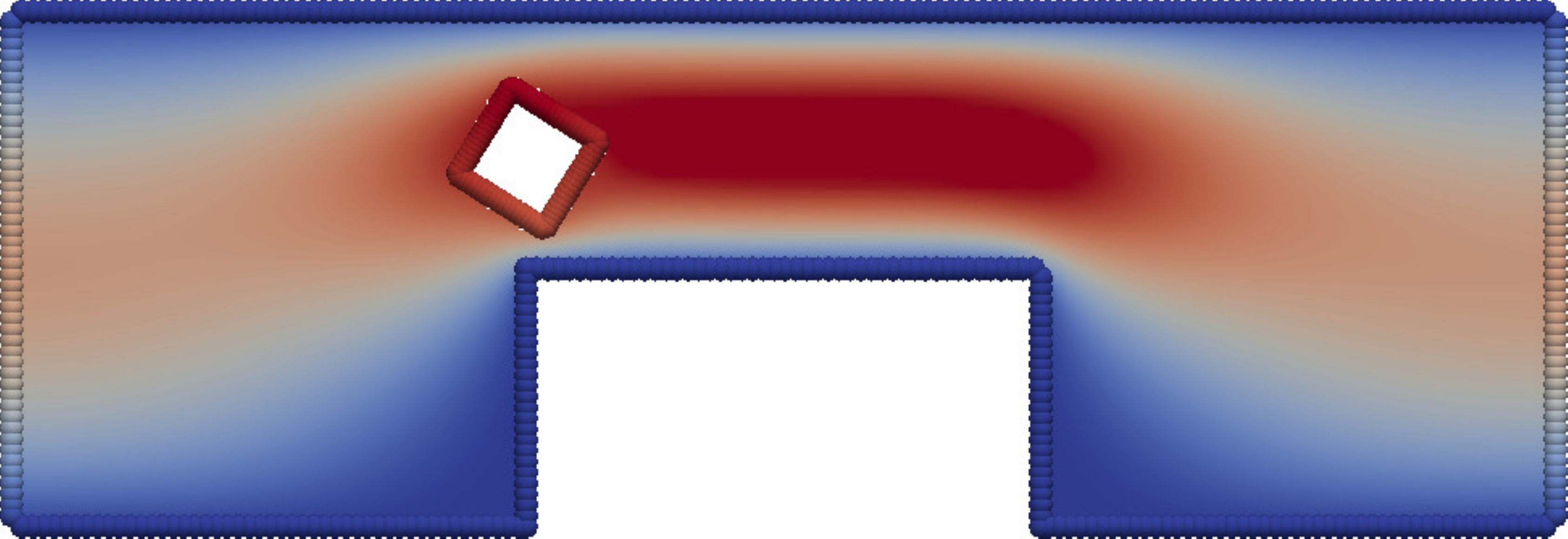}
  \includegraphics[width=0.6\linewidth]{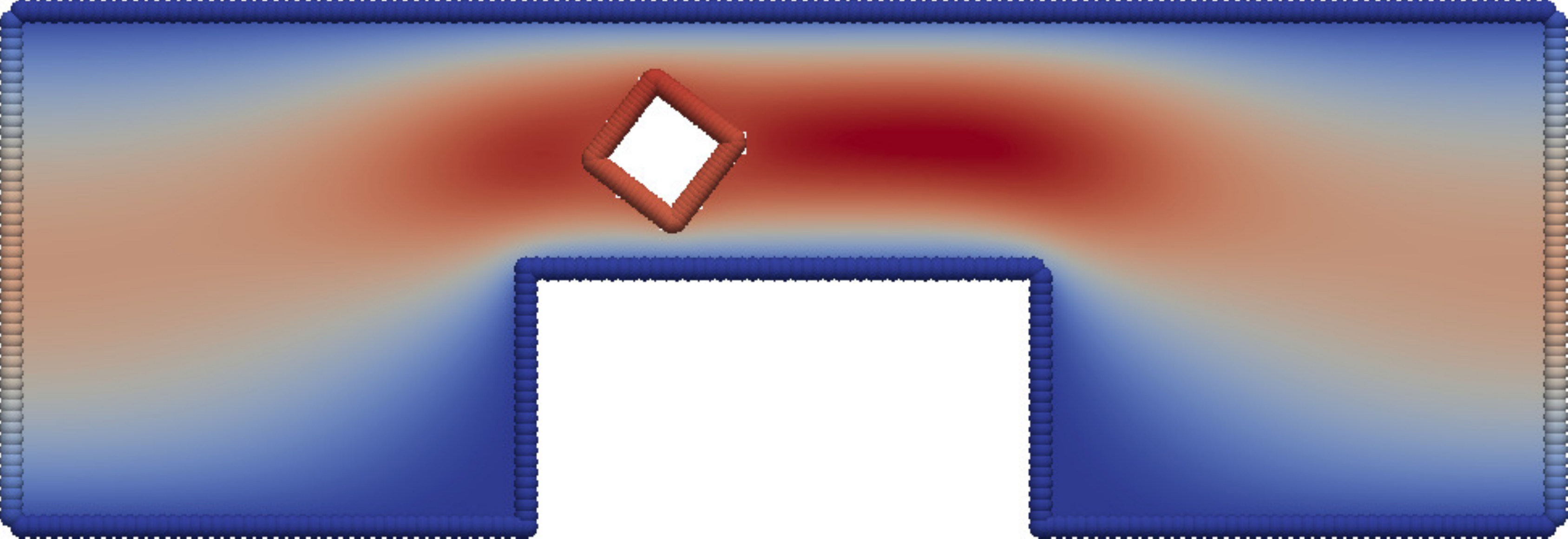}
  \includegraphics[width=0.6\linewidth]{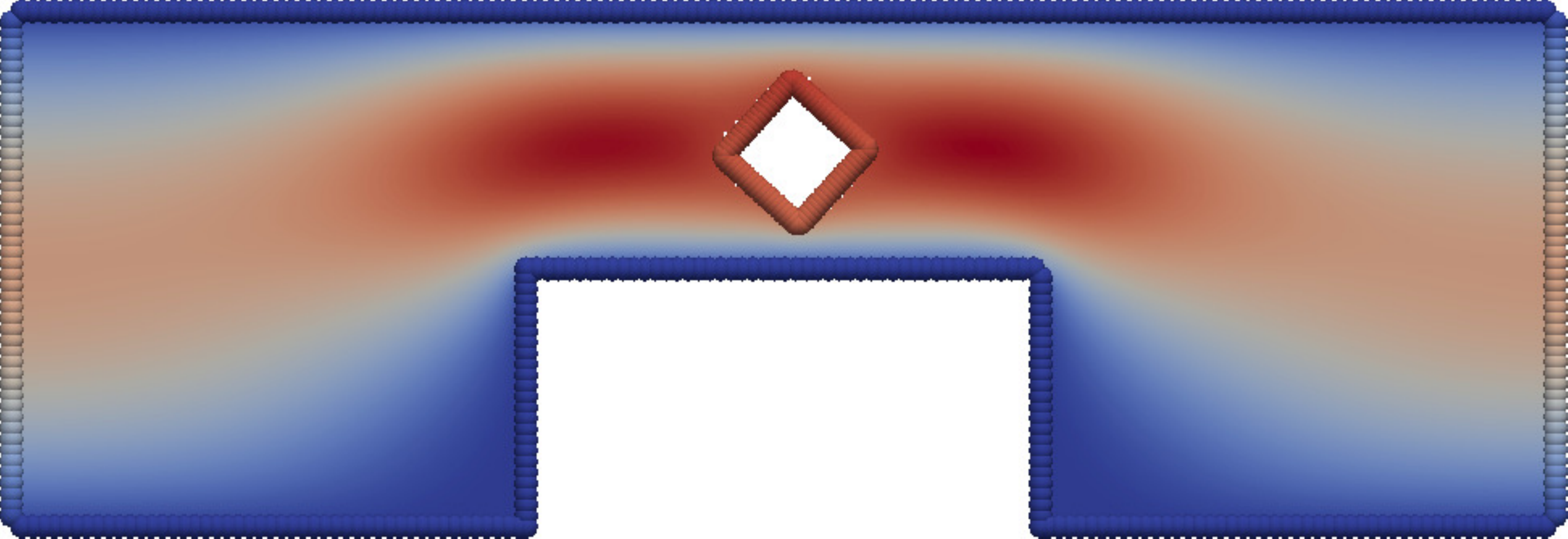}
  \caption{Trajectory of a square colloid passing through a notched microchannel.}
  \label{fig:notched}
\end{figure}

With the Stokes solver and its coupling to colloid dynamics validated, we finally have a means of evaluating the right hand side of Equation \ref{colloidEvolution2} so that the hydrodynamic forces can be integrated numerically to calculate colloid trajectories. We consider the case where two cylindrical particles of radius $a$ are placed in a channel of dimension $L\times H$ with positions $\mathbf{X}_1 = (-L_h,L_v)$ and $\mathbf{X}_2 = (L_h,-L_v)$ and flow is driven by a shear rate $\gamma_0$ (Figure \ref{shearFlow}). For the case where $\frac{a}{L+H} \rightarrow 0$, the trajectory of the particles is governed by a non-linear ODE that can be integrated numerically to provide a highly accurate benchmark solution \cite{arp1977kinetics,bian2014splitting}. This case has been used to validate lower-order methods for simulating colloidal suspensions in the past, however in this case we again stress that no lubrication models are used here. An accurate solution will be obtained directly by solving the Stokes equations using adaptivity to resolve the narrow lubrication gap between particles. We consider $a=1$,$L_h=1.5$ and the family of trajectories with initial condition given by $L_v \in \left\{2,1,0.5,0.25\right\}$ and also the configuration $\mathbf{X}_1 = (-1.2, 0)$,$\mathbf{X}_2 = (1.2, 0)$. The exact solution for this problem predicts that the first set of initial configurations will give rise to open trajectories, while the second will lead to closed orbits in which the two particles will circle each other indefinitely. The shear rate $\dot{\gamma}$ is imposed by enforcing a Couette flow on $\Omega_d$
\begin{equation}
\mathbf{w} = \dot{\gamma} y
\end{equation}
To reproduce the exact solution, it is necessary to take a sufficiently large domain that wall effects do not impact the particle trajectories; we have selected $L=H=40a$. From experiments and the reference solution however, it can be shown that, for the given choice of colloid configuarations, as the particles pass each other the lubrication gap can be as small as $a/50$. That presents a ratio of relevant length scales of $\sim 2000$, and following the previously defined heuristic of five particles per lengthscale would lead to a simulation with $10^{10}$ particles if a uniform discretization lengthscale were used. We will address this by using relatively aggresive adaptivity, using a $\{32^{-1},M_r,5\}$ family of fourth-order discretizations and choose $M_r$ so that there are approximately five points across the lubrication gap. This leads to simulations of roughly ten thousand particles with very fine resolution in the gap, however because the monolithic Stokes/colloid system is completely implicit no stability condition is imposed, and $\Delta t = 0.1$ provides an accurate resolution of the particle trajectories when compared to the benchmark solution (Figure \ref{trajectory}). One of the key benefits of this meshless approach is that a colloidal particle of arbitrary shape can easily be handled. We highlight this capability in Figure \ref{fig:test}, where we demonstrate results for the same problem but using square colloids of unit size. To demonstrate the capability to handle arbitrary domains, we present results in which a square colloid passes through a notched microchannel in Figure \ref{fig:notched}.

\begin{figure}
  \includegraphics[width=.9\linewidth]{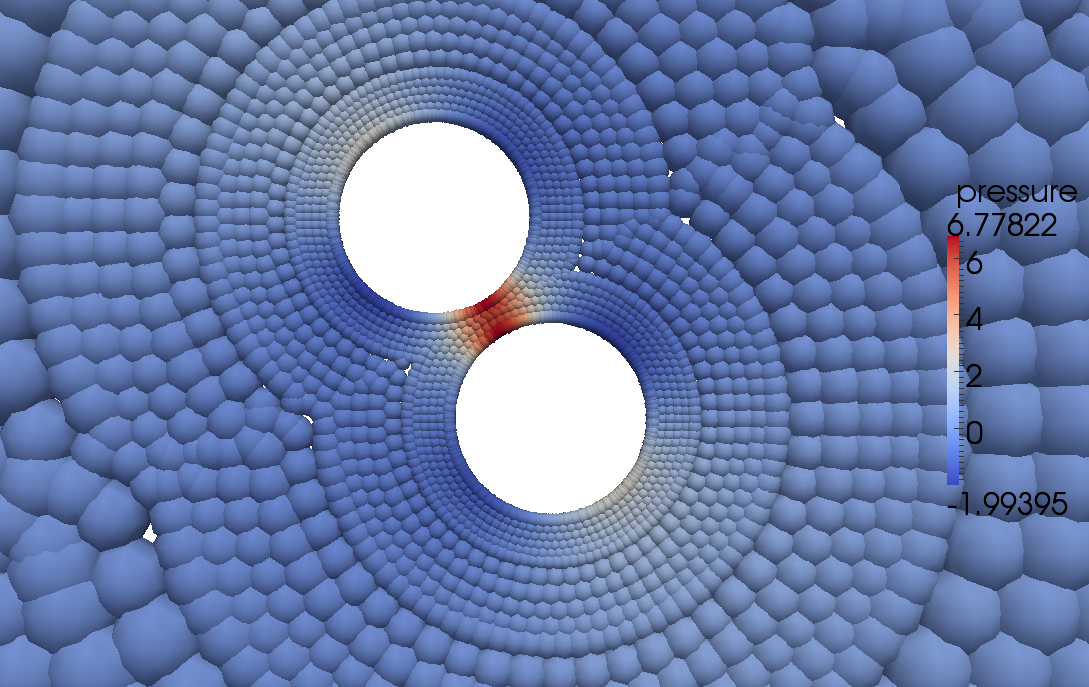}
  \includegraphics[width=.9\linewidth]{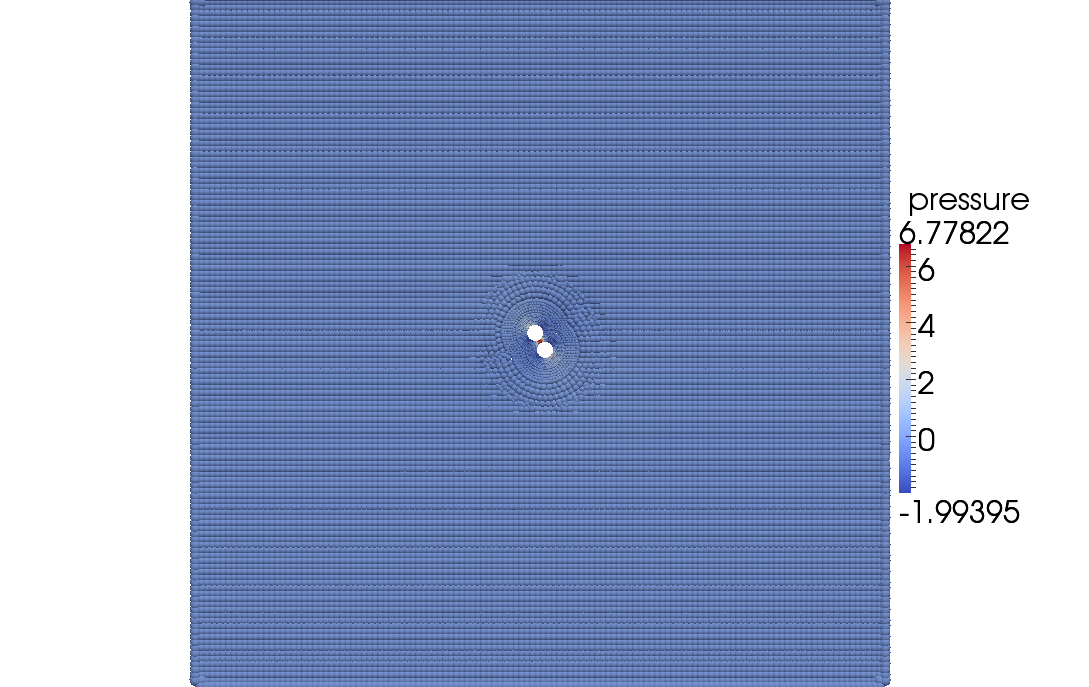}
  \caption{Point adaptivity for colloids interacting under shear flow. Particle adaptivity visualized by rendering spheres at each point with diameter proportional to $\epsilon_i$. First few layers of refinement in vicinity of lubrication layer (top), and complete extant of domain (bottom). We refer to online article for color scale.}
  \label{shearFlow2}
\end{figure}


\section{Conclusion}

In this work we have presented a new compatible meshless discretization for the steady Stokes equations. By combining a staggered MLS discretization for the pressure together with a divergence-free approximation for the velocity, we obtained a stable strong-form meshless discretization that we have shown is able to be efficiently preconditioned using an approximate Schur complement approach. A numerical comparison to a benchmark flow for which analytic solutions are available has shown that the approach is able to achieve equal high-order convergence for both the velocity and the pressure. Equal-order convergence of both pressure and velocity is a property that generally is only achieved by compatible finite element methods. To our knowledge, this marks the first time a meshless method has been able to achieve such high-order convergence while maintaining a sparse discretization and $O(N)$ computational complexity. After demonstrating the accuracy and efficiency of the approach, we used analytic solutions to systematically benchmarked the necessary components to use the scheme to study colloidal suspension flows: evaluation of lubrication forces at the colloid surface, calculation of colloid degrees of freedom under equilibrium conditions, and integration of the particle velocity to obtain colloid trajectories under the flow.

We have demonstrated through these benchmarks that the new discretization provides a new framework for studying suspension flows. Our new approach maintains the flexibility of immersed boundary and low-order meshless methods, while obtaining the high-order of accuracy and sharp treatment of boundary effects that previously could only be obtained via ALE formulations of high-order mesh-based discretizations. Unlike these mesh-based methods however, our meshless approach does not require computationally expensive mesh maintenance at each timestep. The cost of generating a completely new point cloud at every timestep is negligible when compared to the cost of solving the linear system.

A key feature of the results in this work is the ability to accurately resolve lubrication effects for arbitrary geometry directly as a solution to the Stokes equations. By introducing adaptivity in the narrow gap between particles, lubrication forces can be directly resolved without introducing artificial repulsive potentials between colloids. In practice, it is unlikely that such an approach would be ideal - such repulsive potentials have been shown to be an effective tool for reducing computational cost, particularly for physical problems where surface roughness plays an important role and it is unrealistic to pursue such a large degree of adaptivity. In the current work we stress this first principles capability to highlight a case where high-order methods and the a sharp interface representation allow us to approach an inherently multiscale problem.

Having demonstrated the capabilities and robustness of the discretization, we are currently developing a massively parallel three-dimensional implementation of the scheme. The authors' have previously developed a software library coupling the LAMMPS\cite{plimpton1995fast} and Trilinos \cite{heroux2003overview} libraries to enable the development of sparse discretizations resulting from meshless methods\cite{traskCMAME}. This will be adapted to the current scheme to study multiphysics problems in dense suspension flows. We reserve a rigorous benchmarking of the performance of the current scheme for this future work.

While we have focused in this work on the solution of the steady Stokes equations and applications related to suspension flows, we note that the development of consistent projection methods for the Navier-Stokes problem follows from the stable solution of the steady Stokes problem. In a finite element context, the stability of projection schemes has been shown to rely on the inf-sup compatibility between choices of pressure and velocity spaces\cite{liu2007stability}. A number of meshless discretizations have pursued these projection schemes, and we leave a discussion of the impact of our new discretization on such methods for future work.





\bibliography{JCP.bbl}






%

\end{document}